\documentclass[
  acmsmall, % default: `manuscript`; `acmsmall` is used for TOMS;
  screen,
]{acmart}

\newcommand{\coeff}{u}
\newcommand{\coeffb}{v}
\newcommand{\coeffc}{w}
\newcommand{\order}{{\mathcal O}}
\newcommand{\R}{\mathbb{R}}

\newcommand{\revA}[1]{{#1}}
\newcommand{\revB}[1]{{#1}}
\newcommand{\revC}[1]{{#1}}

\usepackage{xspace}
\newcommand{\bseries}{{BSeries.jl}\xspace}
\newcommand{\pybs}{\texttt{pybs}\xspace}

\usepackage{xcolor}
\usepackage{subfig}
\usepackage{pdflscape}
\usepackage{algorithm}
\usepackage{algpseudocode}

\usepackage{forest}
\forestset{
  */.style={
    delay+={append={[]},}
  },
  rooted tree/.style={
    for tree={
      grow'=90,
      parent anchor=center,
      child anchor=center,
      s sep=2.5pt,
      if level=0{
        baseline
      }{},
      delay={
        if content={*}{
          content=,
          append={[]}
        }{}
      }
    },
    before typesetting nodes={
      for tree={
        circle,
        fill,
        minimum width=3pt,
        inner sep=0pt,
        child anchor=center,
      },
    },
    before computing xy={
      for tree={
        l=5pt,
      }
    }
  },
  big rooted tree/.style={
    for tree={
      grow'=90,
      parent anchor=center,
      child anchor=center,
      s sep=2.5pt,
      if level=0{
        baseline
      }{},
      delay={
        if content={*}{
          content=,
          append={[]}
        }{}
      }
    },
    edge={
        semithick,
        -Latex
    },
    before typesetting nodes={
      for tree={
        circle,
        fill,
        minimum width=5pt,
        inner sep=0pt,
        child anchor=center,
      },
    },
    before computing xy={
      for tree={
        l=15pt,
      }
    }
  },
  my node/.style={
    circle, fill=#1,
  },
  gray node/.style={
    my node=gray, inner sep=#1mm,
  },
  red node/.style={
    my node=red, inner sep=#1mm,
  },
}
\DeclareDocumentCommand\rootedtree{o}{\Forest{rooted tree [#1]}}
\DeclareDocumentCommand\bigrootedtree{o}{\Forest{big rooted tree [#1]}}

\usepackage{alphabeta}

\usepackage{siunitx}
\begin{document}

\title{Computing with B-series}

\author{David I. Ketcheson}
\email{david.ketcheson@kaust.edu.sa}
\affiliation{%
  \institution{CEMSE Division, King Abdullah University of Science \& Technology (KAUST)}
  \streetaddress{4700 KAUST}
  \city{Thuwal}
  \state{Makkah}
  \country{Saudi Arabia}
  \postcode{23955}
}
\author{Hendrik Ranocha}
\email{mail@ranocha.de}
\affiliation{%
  \institution{Applied Mathematics, University of M\"unster}
  \streetaddress{Orl\'eans-Ring 10}
  \postcode{48149}
  \city{M\"unster}
  \country{Germany}
}
\authornote{Current address: Applied Mathematics, University of Hamburg, Bundesstr. 55, 20146 Hamburg, Germany}

\begin{abstract}
We present BSeries.jl, a Julia package for the computation and manipulation of B-series, which are a
versatile theoretical tool for understanding and designing discretizations of
differential equations.
We give a short introduction to the theory of B-series
and associated concepts and provide examples of their use, including method composition
and backward error analysis.  The associated software is highly performant
and makes it possible to work with B-series of high order.
\end{abstract}




\maketitle

\section{Introduction}
B-series are a theoretical tool related to rooted trees and developed
originally for the analysis of Runge-Kutta methods \cite{butcher1972algebraic,hairer1974butcher}.
These topics are described beautifully in sources such as
\cite{hairer1993,butcher2003,butcher2009b,hairer2013geometric,Chartier2010,butcher2021b}, and we follow
the notation and presentation from those sources as much as possible.
The purpose of the present work is to provide an introduction to
this theory in the context of a computational software package that
facilitates both understanding and using B-series.

There exist a number of computational tools for working with Runge-Kutta order conditions
\cite{sofroniou1994symbolic,gruntz1995symbolic,bornemann2001runge,bornemann2002runge,famelis2004symbolic,cameron2006matlab,ketcheson2020nodepy,knothbseries}; these use B-series implicitly or incidentally for this specific goal.
In contrast, \bseries was written for the purpose of computing and manipulating
B-series themselves, facilitating analysis that could not be done with the above packages.
The only similar library of which we are aware is the Python code {\pybs } \cite{sundklakk2015pybs}.
We developed \bseries after also testing an initial prototype in Python \cite{ketcheson2021bseries}.
\bseries is written to be highly performant in order to enable
working with high-order B-series, and includes a more comprehensive set of
features such as the ability to generate elementary differentials for specific
equations and thus form the full B-series for a given integrator and ODE, as well
as the generation of modifying integrators \revA{(see Section \ref{sec:modifying-integrators})}.

It is possible to interface with \bseries at a range of levels.  Some of the features
implemented in \bseries, starting with lower-level functionality and working up, are:
\begin{itemize}
    \item A realization of B-series as maps from rooted trees to the real numbers
    \item \revA{Generalized B-series based on colored rooted trees}
    \item \revA{Vector space operations on B-series}
    \item Computation of the B-series of any Runge-Kutta method
    \item \revA{Computation of the B-series of any additive Runge-Kutta method}
    \item \revA{Computation of the B-series of any Rosenbrock-Wanner (ROW) method}
    \item \revA{The composition and substitution laws for B-series, based on ordered subtrees and partitions of trees}
    \item \revA{Computation of the modified equation for B-series methods, or for a combination
            of method and ODE}
    \item \revA{Computation of a modifying integrator for B-series methods, or for a combination
            of method and ODE}
    \item Human-readable output of computed B-series
\end{itemize}
All operations can be performed using symbolic inputs, for instance to study parameterized
families of methods. Currently, \bseries focuses on working with B-series
truncated at a fixed but arbitrary order.

\bseries \cite{ranocha2021bseries} is written in Julia \cite{bezanson2017julia}
and builds upon the Julia package RootedTrees.jl \cite{ranocha2019rootedtrees},
which we initiated earlier and enhanced while creating \bseries.
In addition to features to enable the functionality of \bseries listed above,
RootedTrees.jl comes with the following features.
\begin{itemize}
  \item Computation of functions on \revA{(colored)} rooted trees
  \item \revA{(Additive)} Runge-Kutta order conditions
  \item Human-readable output of rooted trees in various text and graphics formats
\end{itemize}
\bseries supports different packages for symbolic computations, including
SymPy \cite{meurer2017sympy} via SymPy.jl\footnote{\url{https://github.com/JuliaPy/SymPy.jl}},
Symbolics.jl\footnote{\url{https://github.com/JuliaSymbolics/Symbolics.jl}} \cite{gowda2021high},
and SymEngine\footnote{\url{https://github.com/symengine/symengine}} via SymEngine.jl\footnote{\url{https://github.com/symengine/SymEngine.jl}}.

Code\footnote{\url{https://github.com/ranocha/BSeries.jl}} and
documentation\footnote{\url{https://ranocha.github.io/BSeries.jl/stable/}} for
\bseries is available online, including an API reference and
examples.  The package includes a set of tests that currently cover 99\% of the
code and are run automatically via continuous integration.  The code is distributed
under an MIT license and contributions are welcomed.
All code necessary to reproduce the examples shown in this article is available
in our reproducibility repository \cite{ketcheson2021computingRepro}.

\subsection{Rooted trees}
A rooted tree is an acyclic graph with one node designated as the root.
We denote the set of rooted trees with $k$ nodes by $T_k$, and the set
of all rooted trees by $T \revA{= \bigcup_{k \in \mathbb{N}} T_k}$.  The first few sets $T_k$ are
\begin{align}
    T_1 & = \{ \rootedtree[] \} \\
    T_2 & = \{ \rootedtree[[]] \} \\
    T_3 & = \{ \rootedtree[[],[]], \rootedtree[[[]]] \} \\
    T_4 & = \{ \rootedtree[[],[],[]], \rootedtree[[[],[]]], \rootedtree[[],[[]],], \rootedtree[[[[]]]] \}
\end{align}
Note that here we are concerned with unlabeled trees, so we do not distinguish between
$\rootedtree[[[]],[]]$ and $\rootedtree[[],[[]]]$.  We will sometimes write a tree
in terms of its children:
$$
    \tau = [ \tau_1, \dots, \tau_m ].
$$
Here $\tau_1, \dots, \tau_m$ are the trees remaining when we remove
the root of $\tau$; we will sometimes refer to them as the \emph{subtrees} of
$\tau$.  We will make use of the following functions on trees:
\begin{itemize}
\item As discussed already, the \emph{order} of a tree, denoted by $|\tau|$, is the number of nodes it contains.
\item The \emph{symmetry} of a tree is given recursively by
\begin{align}
    \sigma(\rootedtree[])=1, \quad \quad \quad \sigma(\tau) = \sigma(\tau_1)\cdots \sigma(\tau_m) \cdot \mu_1!\mu_2!\cdots
\end{align}
where the integers $\mu_1,\mu_2,\dots$ are the numbers of identical trees
among $\tau_1,\dots,\tau_m$.  Thus for instance $\sigma(\rootedtree[[[]]])=1$
and $\sigma(\rootedtree[[],[]])=2$.
\item The \emph{density} of a tree is defined recursively by
\begin{align*}
    \gamma(\rootedtree[])& =1,   & \gamma(\tau) & = |\tau|\gamma(\tau_1)\cdots\gamma(\tau_m).
\end{align*}
\end{itemize}

\begin{table}
  \caption{The first few rooted trees and some of their properties.}
  \label{tbl:trees}
  \begin{tabular}{cccclll}
    \toprule
    $\tau$ & $|\tau|$ & $\sigma(\tau)$ & $\gamma(\tau)$ & $F(\tau)$ & $\Phi_j(t)$ & Level sequence \\
    \midrule
$\rootedtree[]$ & 1 & 1 & 1 & $f^j$ & $1$ & $[1]$ \\
$\rootedtree[[]]$ & 2 & 1 & 2 & $f^j_{k} f^k$ & $a_{ jk }$ & $[1, 2]$ \\
$\rootedtree[[], []]$ & 3 & 2 & 3 & $f^j_{kl} f^k f^l$ & $a_{ jk } a_{ jl }$ & $[1, 2, 2]$ \\
$\rootedtree[[[]]]$ & 3 & 1 & 6 & $f^j_{k} f^k_{l} f^l$ & $a_{ jk } a_{ kl }$ & $[1, 2, 3]$ \\
$\rootedtree[[], [], []]$ & 4 & 6 & 4 & $f^j_{klm} f^k f^l f^m$ & $a_{ jk } a_{ jl } a_{ jm }$ & $[1, 2, 2, 2]$ \\
$\rootedtree[[], [[]]]$ & 4 & 1 & 8 & $f^j_{km} f^k_{l} f^l f^m$ & $a_{ jk } a_{ jm } a_{ kl }$ & $[1, 2, 3, 2]$ \\
$\rootedtree[[[], []]]$ & 4 & 2 & 12 & $f^j_{k} f^k_{lm} f^l f^m$ & $a_{ jk } a_{ kl } a_{ km }$ & $[1, 2, 3, 3]$ \\
$\rootedtree[[[[]]]]$ & 4 & 1 & 24 & $f^j_{k} f^k_{l} f^l_{m} f^m$ & $a_{ jk } a_{ kl } a_{ lm }$ & $[1, 2, 3, 4]$ \\
  \bottomrule
\end{tabular}
\end{table}

\revA{In this paper we represent trees graphically.  In code, we represent them via level sequences.

\begin{definition}
  The \emph{level} of a vertex $v$ of a rooted tree $t$ is unity if $v$ is the
  root of $t$. Otherwise, it is one greater than the level of the parent of $v$.
  A \emph{level sequence} $(l_1, l_2, \dots, l_n)$ of a rooted tree $t$ with $n$
  vertices is obtained by recording the level of each vertex of $t$ in a
  depth-first search starting at the root of $t$.
\end{definition}

The set of rooted trees of a given order can be generated using
the constant time algorithm of \cite{beyer1980constant}.
Ordered rooted trees can be compared via lexicographical ordering of their level
sequences. However, the level sequence representation of an unordered rooted tree
is not unique. For many algorithms, including the generation of rooted trees
\cite{beyer1980constant}, it is advantageous to use the \emph{canonical representation}
of a rooted tree $t$ based on the lexicographically biggest level sequence associated
to $t$.
Level sequences for a few trees are given in Table \ref{tbl:trees}.}

\subsection{Discretizations of ordinary differential equations}
The study of B-series arose from the analysis of numerical solutions of first-order
systems of differential equations
\begin{align} \label{ode}
    \dot{y} = f(y),
\end{align}
where $y=y(t)$ and $\dot{y}$ denotes the derivative of $y$ with respect to $t$.
We are focused on the initial value problem where $y(t=0)$ is specified, but
we will usually suppress $t$ to shorten the  expressions appearing later in
this work.

The solution of \eqref{ode} is often approximated by a Runge-Kutta method,
which produces a sequence of approximations $y_n \approx y(nh)$
\begin{subequations} \label{RK}
\begin{align}
    Y_i & = y_n + h \sum_{j=1}^s a_{ij} f(Y_j), \\
    y_{n+1} & = y_n + h \sum_{j=1}^s b_j f(Y_j),
\end{align}
\end{subequations}
where $h$ is the step size and $A, b$ are the coefficients that
specify the method.  \revA{For non-autonomous problems (with $f = f(y,t)$),
a Runge-Kutta method requires the additional coefficients $c_1, \dots, c_s$
and the derivatives on the right-hand side of \eqref{RK} are given by
$f(Y_j,t_n+c_j h)$.
We assume throughout this work that $c_i = \sum_j a_{ij}$.  For non-autonomous problems,
this guarantees that the method gives the same result whether it is written in
autonomous or non-autonomous form.  B-series analysis is most conveniently
applied to the autonomous form.}

\revC{A Runge-Kutta method is often represented by the Butcher tableau
\begin{align}
\begin{array}{c|c}
    c & A \\ \hline
      & b^T
\end{array}
\end{align}
where $b, c \in \R^s$ and $A\in \R^{s\times s}$.}

Direct analysis of the numerical solution based on \eqref{RK}
is challenging, since it involves nested evaluations of $f$.  Instead,
one can write the solution as an (infinite) power series in $h$ and involving
derivatives of $f$.
This is known as a B-series, and takes the form
\begin{align} \label{bseries}
    B(\coeff, hf, y) & := \coeff(\varnothing) y + \sum_{\tau \in T} \frac{h^{|\tau|}}{\sigma(\tau)} \coeff(\tau) F_f(\tau)(y),
\end{align}
where $T$ is the set of all rooted trees and $|\tau|$ is the \emph{order} of
the rooted tree $\tau$, i.e., the number of nodes in $\tau$.
The factor $F_f(\tau)(y)$ involves the derivatives of $f$ and is referred to
as an \emph{elementary differential}.

We will work with B-series that represent a map:
\begin{align} \label{Bmap}
    y_{n+1} & = B(\coeff, hf, \revC{y_n})
\end{align}
and B-series that represent a flow:
$$
    \dot{y}(t) = B(\coeff, hf, y).
$$
In the former case consistency demands that $\coeff(\varnothing) = 1$,
while in the latter we must have $\coeff(\varnothing) = 0$.

\subsection{Elementary differentials}
The derivatives of $f$ with respect to $t$
involve various products of the tensors of partial derivatives of $f$.
The derivatives of order $k$ are in one-to-one correspondence with the set of
rooted trees with $k$ nodes.
The factor $F_f(\tau)(y)$ in \eqref{bseries} is an \emph{elementary differential}.
\revA{Here we primarily follow the notation established in \cite{hairer1993,hairer2013geometric}, which the reader may
wish to consult for more details.}
For simplicity we will usually omit the subscript $f$.
The elementary differentials for a given function $f$ are defined recursively by
\begin{align}
    F(\rootedtree[])(y) & = f(y) \\
    F(\tau)(y) & = f^{(m)}\left(F(\tau_1)(y), \cdots, F(\tau_m)(y)\right) & \text{for } \tau=[\tau_1,\dots,\tau_m].
\end{align}
From a programming perspective it is perhaps clearer to write the
elementary differentials in Einstein summation form; e.g.
$$
    F\Bigl(\rootedtree[[[]],[]]\Bigr)(y)
    = f''(f' f, f)
    = f^j_{kl} f^l_m f^m f^k.
$$
Here subscripts denote differentiation and repeated subscripts imply summation.

\section{Applications of B-series to Runge-Kutta methods}

\subsection{Expression of the numerical solution map}
As mentioned above, B-series can be used to express the operation of a numerical
method in terms of an infinite series.
The B-series of a Runge-Kutta
method with coefficients $A, b$ is given by
\begin{align} \label{RKmap}
    y_{n+1} = B(\Phi,hf,y_n),
\end{align}
where the \emph{elementary weights} $\Phi(\tau)$ are defined via
\revC{
\begin{align}
    \Phi_j(\rootedtree[]) & = 1 \\
    \Phi_j(\tau) & = \sum_{k_1,\dots,k_m}^s (a_{j{k_1}} \Phi_{k_1}(\tau_1)) \cdots (a_{j{k_m}} \Phi_{k_m}(\tau_m)) & \text{for } \tau=[\tau_1,\dots,\tau_m] \\
    %\Phi_j(\tau) & = \sum_{k=1}^s a_{jk} \prod_{l=1}^m \Phi_k(\tau_l) \\
    \Phi(\tau) & = \sum_{j=1}^s b_j \Phi_j(\tau),
\end{align}
where $s$ is the number of stages in the method and $m$ is the number of subtrees
of $\tau$.}
Some examples are given in Table \ref{tbl:trees}.

\subsubsection{Example: \revA{explicit} second-order two-stage Runge-Kutta method}
Using the expressions above, we can write down the B-series of a given
method or even a whole family of methods.  For instance, the
one-parameter family of \revA{explicit} two-stage, second-order Runge-Kutta methods
is given by the Butcher arrays
\begin{align} \label{rk22}
\begin{array}{c|cc}
 0  & 0 & \\
\frac{1}{2\alpha} & \frac{1}{2\alpha} & 0 \\[0.1cm]
\hline
 & 1 - \alpha & \alpha\\
\end{array}
\end{align}
When this method is applied to a system of ODEs, the leading terms in the resulting
solution map \eqref{Bmap} are
\begin{multline}
  y_{n+1}
  =
  y_n
  + h F_{f}\mathopen{}\left( \rootedtree[] \right)\mathclose{}(y_n)
  + \frac{1}{2} h^{2} F_{f}\mathopen{}\left( \rootedtree[[]] \right)\mathclose{}(y_n)
  + \frac{1}{8 \alpha} h^{3} F_{f}\mathopen{}\left( \rootedtree[[][]] \right)\mathclose{}(y_n)
  \\
  + \frac{1}{48 \alpha^{2}} h^{4} F_{f}\mathopen{}\left( \rootedtree[[][][]] \right)\mathclose{}(y_n)
  + \frac{1}{384 \alpha^{3}} h^{5} F_{f}\mathopen{}\left( \rootedtree[[][][][]] \right)\mathclose{}(y_n) + \order(h^6)
\end{multline}
which can be obtained via
\begin{verbatim}
    using BSeries, SymPy, Latexify
    α = symbols("α", real=true)
    A = [0 0; 1/(2*α) 0]; b = [1-α, α]; c = [0, 1/(2*α)]
    coefficients = bseries(A, b, c, 5)
    println(latexify(coefficients, cdot=false))
\end{verbatim}
Using such an expression, one could study for instance what choices of $\alpha$ will
suppress (to leading orders) certain components of the error that are related to
particular elementary differentials.
Determining expressions of this kind by hand would be extremely tedious and error-prone.

\subsection{Order conditions}
The exact solution map for the flow $\dot{y}(t)=f(y)$ is
\begin{align}
    y(t+h) = B(e,hf,y(t))
\end{align}
where the B-series coefficients of the exact solution are \cite[Thm. II.2.6]{hairer1993}
\begin{align} \label{exact}
    e(\tau) = \frac{1}{\gamma(\tau)}.
\end{align}
To determine the accuracy of a Runge-Kutta method, we should compare its
B-series coefficients to the values of $e(\tau)$.
Comparing \eqref{RKmap} with \eqref{exact} yields the conditions for a Runge-Kutta method
to have order $p$:
$$
    \Phi(\tau) = e(\tau) \quad \forall \tau \in T \text{ with } |\tau|\le p.
$$

\subsection{Composition of methods}
The composition of two B-series integrators is another B-series integrator.
This can be used to understand the behavior of a numerical solution when
different methods are applied in sequence, or to design new numerical methods.
The composition operation is the basis for viewing the Runge-Kutta methods as
a group, as introduced originally by Butcher \revA{\cite{butcher1972algebraic}}.  The algebraic relation
between the B-series of the component methods and their composition is
\begin{align} \label{compshort}
    B(\coeff, hf, B(\coeffb, hf, y)) = B(\coeffb\cdot \coeff, hf, y)
\end{align}
where the composition operation is
\begin{align} \label{composition}
    (\coeffb \cdot \coeff)(\varnothing) & = \coeff(\varnothing) \\
    (\coeffb \cdot \coeff)(\tau) & = \sum_{s\in S(\tau)} \coeffb(\tau \setminus s) \coeff(s_\tau).
\end{align}

\revA{\begin{definition} {\bf Ordered subtrees, and the associated forest.}
An ordered subtree of a rooted tree $\tau$ is a subset of the nodes of $\tau$ such that the
root is included and the nodes are connected (by edges of $\tau$).
%The number of ordered subtrees is generally less than $2^{|\tau|-1}$ since not all set of nodes leave a connected tree.
Taking these nodes along with the edges in $\tau$
that connect them gives the rooted tree $s_\tau$.  Removing these nodes and the adjacent
edges from $\tau$ gives the forest $\tau \setminus s$.
The set of all ordered subtrees of $\tau$ is denoted by $S(\tau)$.
These concepts are used in the definition of the composition law \eqref{composition}.
Examples are given in Table \ref{tbl:subtrees}.  Note that the empty set is also
allowed as the contents of $s_\tau$ or $\tau \setminus s$.
\end{definition}}

Although the formula \eqref{compshort} is written with equal step sizes for
the two methods, we can use it to compute compositions with different
step sizes by noting that
$$
    B(\coeff, \mu h f, y) = B(\coeffc, hf, y)
$$
where $\coeffc(\tau) = \mu^{|\tau|} \coeff(\tau)$.

\subsubsection{Composition of two second-order methods}
As a first example, consider the problem of composing two 2-stage, 2nd-order
Runge-Kutta methods of the form \eqref{rk22} (with the same step size
but possibly different values of $\alpha$) in order to obtain a method
of higher order.  After generating the B-series for each method and
storing them in \verb|series1| and \verb|series2|, the
composition is computed with the single line of code
\begin{verbatim}
    comp = compose(series1, series2, normalize_stepsize=true)
\end{verbatim}
The initial terms of the resulting series are
$$y + hf(y) + \frac{1}{2} h^{2}
F_{f}\mathopen{}\left( \rootedtree[[]] \right)\mathclose{}(y) + \frac{1}{8} h^{3}
F_{f}\mathopen{}\left( \rootedtree[[[]]] \right)\mathclose{}(y) + \left(\frac{1}{8} +
\frac{1}{64 \alpha_1} + \frac{1}{64 \alpha_2}\right) h^{3} F_{f}\mathopen{}\left(
\rootedtree[[][]] \right)\mathclose{}(y),$$
from which we see that the term corresponding to \rootedtree[[[]]] is independent
of the values $\alpha_1, \alpha_2$ and different from $1 / \gamma(t) = 1 / 6$,
the corresponding coefficient of the exact solution B-series. Thus, any such
composition will have order equal to two.

\revB{
\subsubsection{Effective order}
The concept of effective order is based on the idea of composing RK methods
\cite{butcher1969effective}.
Here, we present the example of Butcher's method of effective order 5. This is
a fourth-order Runge-Kutta method that results in a fifth-order method when
composed with a special starting and finishing procedure.

The Butcher tableau for the main method is:
\begin{align}
\begin{array}{c|ccccc}
 &  &  &  &  & \\
\frac{1}{5} & \frac{1}{5} &  &  &  & \\
\frac{2}{5} &  & \frac{2}{5} &  &  & \\
\frac{1}{2} & \frac{3}{16} &  & \frac{5}{16} &  & \\
1 & \frac{1}{4} &  & - \frac{5}{4} & 2 & \\
\hline
 & \frac{1}{6} &  &  & \frac{2}{3} & \frac{1}{6}\\
\end{array}
\end{align}
Those of the starting and finishing method are:
\begin{align}
\begin{array}{c|ccccc}
 &  &  &  &  & \\
\frac{1}{5} & \frac{1}{5} &  &  &  & \\
\frac{2}{5} &  & \frac{2}{5} &  &  & \\
\frac{3}{4} & \frac{75}{64} & - \frac{9}{4} & \frac{117}{64} &  & \\
1 & - \frac{37}{36} & \frac{7}{3} & - \frac{3}{4} & \frac{4}{9} & \\
\hline
 & \frac{19}{144} &  & \frac{25}{48} & \frac{2}{9} & \frac{1}{8}\\
\end{array}
& &
\begin{array}{c|ccccc}
 &  &  &  &  & \\
\frac{1}{5} & \frac{1}{5} &  &  &  & \\
\frac{2}{5} &  & \frac{2}{5} &  &  & \\
\frac{3}{4} & \frac{161}{192} & - \frac{19}{12} & \frac{287}{192} &  & \\
1 & - \frac{27}{28} & \frac{19}{7} & - \frac{291}{196} & \frac{36}{49} & \\
\hline
 & \frac{7}{48} &  & \frac{475}{1008} & \frac{2}{7} & \frac{7}{72}\\
\end{array}
\end{align}
We can generate the B-series for the main method and check its accuracy as follows:}
\begin{verbatim}
using BSeries

A = [0 0 0 0 0;
     1//5 0 0 0 0;
     0 2//5 0 0 0;
     3//16 0 5//16 0 0;
     1//4 0 -5//4 2 0]
b = [1 // 6, 0, 0, 2 // 3, 1 // 6]
rk_a = RungeKuttaMethod(A, b)
series_a = bseries(rk_a, 6)
order_of_accuracy(series_a)
\end{verbatim}
\revB{
This returns 4 since the order conditions are satisfied for all trees up to order 4, 
but not for trees of order 5, so the method is fourth-order accurate.  
A similar check shows that the starting and finishing methods are third-order accurate.

We generate their composition and compare it with the B-series
of the exact solution as follows:}
\begin{verbatim}
series_comp = compose(series_b, series_a, series_c, normalize_stepsize = true)
order_of_accuracy(series_comp)
\end{verbatim}
\revB{
This returns 5, confirming the fifth-order accuracy of the composed method.
}

\subsection{Backward error analysis}

Backward error analysis is a fundamental tool of numerical analysis.
Instead of studying the difference between the computed solution and
the exact solution (i.e., the \emph{forward error}), one assumes that
the computed solution is the exact solution of a perturbed problem,
and studies how that problem differs from the original one.
This perspective is widely used in numerical linear algebra and
numerical methods for linear partial differential equations (PDEs).
In the latter setting it is referred to as \emph{modified equation analysis}
since the result is a modification of the original PDE (see e.g. \cite{rjl:fdmbook,karam2020pymodpde}).
This approach is less well-known in the context of nonlinear differential
equations, probably because it is much more difficult (or at least tedious)
to apply.  It has mostly been used in the study of symplectic or reversible
methods \cite{hairer2013geometric}.  Our original motivation for developing
\bseries \cite{ranocha2021bseries} was to automate and facilitate
backward error analysis for the numerical solution of nonlinear ODEs.

Here we follow the exposition from \cite{Chartier2010}.
Let $\phi_h$ denote the numerical solution map; i.e.
$$
y_{n+1} = \phi_h(y_n).
$$
We suppose there exists a modified differential equation $\dot{y}=f_h(y)$
such that $\phi_h$ is the exact solution map.  Then if $\phi_h$ can be
written as a B-series $B(\coeff,hf,y)$, the modified differential equation
right-hand side also has a B-series representation.  These two
B-series are related by the fact that substitution of the second one
into the exact solution B-series yields the first one:
$$
    B(e,hf_h,y) = B(e, B(v,hf,\cdot), y) = B(u, hf, y).
$$
Here we have used $v$ to denote the coefficients of the second B-series; notice that
the right-hand side of the modified equation is given by $f_h = h^{-1} B(\coeffb,hf,y)$.  
The foregoing motivates the definition of the substitution operator $\star$:
$$
    B(\coeffb \star e, hf, y) := B(e, B(\coeffb, hf, \cdot), y) = B(\coeff, hf, y).
$$
In order to find the modified equation $\dot{y}=f_h(y)$ we need only to
compute the substitution coefficients
$\coeffb \star e$.  This relies on the skeleton $p_\tau$ and the forest
$\tau \setminus p$ associated with the partitions $p\in P$
of \revA{a tree} $\tau$, defined as follows.

\revA{\begin{definition}
{\bf Partitions of a tree, and the associated forest and skeleton.}
A partition $p$ of a tree $\tau$ is any subset of the edges of the tree.  Since tree
$\tau$ has $|\tau|-1$ edges, it has $2^{|\tau|-1}$ partitions.  The set $\tau \setminus p$
is the forest that is left when the edges of $p$ are removed from $\tau$.  The skeleton
$p_\tau$ is the tree that remains when each tree of $\tau\setminus p$ is contracted
to a single node and the edges of $p$ are put back.
Examples are given in Table \ref{tbl:partitions}.
\end{definition}}

In terms of the above entities, we have (see \cite[Theorem 3.2]{Chartier2010})
\begin{equation} \label{substition}
\begin{aligned}
(\coeffb  \star \coeff)(\varnothing) & = \coeff(\varnothing) \\
(\coeffb  \star \coeff)(\tau) & = \sum_{p\in P(\tau)} \coeffb(\tau \setminus p) \coeff(p_\tau).
\end{aligned}
\end{equation}
The system of equations $(\coeffb\star e)(\tau) = \coeff(\tau)$ can be solved recursively
because, for each tree $\tau$, we have $(\coeffb\star e)(\tau) = \coeffb(\tau) +$ (terms depending on
lower-order trees).

\begin{landscape}

\begin{table}
    \centering
  \caption{All partitions of the tree with level sequence $[1,2,3,2,3]$}
  \label{tbl:partitions}
  \begin{tabular}{ccccccccccccccccc}
    \toprule
    $p$ & \bigrootedtree[[[]],[[]]]
        & \bigrootedtree[[.,edge=dashed[]],[[]]]
        & \bigrootedtree[[.,edge=dashed[]],[.,edge=dashed[]]]
        & \bigrootedtree[[.,edge=dashed[]],[[.,edge=dashed]]]
        & \bigrootedtree[[.,edge=dashed[.,edge=dashed]],[[]]]
        & \bigrootedtree[[.,edge=dashed[]],[.,edge=dashed[.,edge=dashed]]]
        & \bigrootedtree[[.,edge=dashed[.,edge=dashed]],[[.,edge=dashed]]]
        & \bigrootedtree[[.,edge=dashed[.,edge=dashed]],[.,edge=dashed[]]]
        & \bigrootedtree[[.,edge=dashed[.,edge=dashed]],[.,edge=dashed[.,edge=dashed]]]
        & \bigrootedtree[[[]],[.,edge=dashed[]]]
        & \bigrootedtree[[[]],[.,edge=dashed[.,edge=dashed]]]
        & \bigrootedtree[[[.,edge=dashed]],[.,edge=dashed[]]]
        & \bigrootedtree[[[.,edge=dashed]],[.,edge=dashed[.,edge=dashed]]]
        & \bigrootedtree[[[.,edge=dashed]],[[]]]
        & \bigrootedtree[[[]],[[.,edge=dashed]]]
        & \bigrootedtree[[[.,edge=dashed]],[[.,edge=dashed]]]  \\[0.5cm]
    $\tau \setminus p$
        & \bigrootedtree[[[]],[[]]]
        & \bigrootedtree[[]] \bigrootedtree[[[]]]
        & \bigrootedtree[] \bigrootedtree[[]] \bigrootedtree[[]]
        & \bigrootedtree[] \bigrootedtree[[]] \bigrootedtree[[]]
        & \bigrootedtree[] \bigrootedtree[] \bigrootedtree[[[]]]
        & \bigrootedtree[] \bigrootedtree[] \bigrootedtree[] \bigrootedtree[[]]
        & \bigrootedtree[] \bigrootedtree[] \bigrootedtree[] \bigrootedtree[[]]
        & \bigrootedtree[] \bigrootedtree[] \bigrootedtree[] \bigrootedtree[[]]
        & \bigrootedtree[] \bigrootedtree[] \bigrootedtree[] \bigrootedtree[] \bigrootedtree[]
        & \bigrootedtree[[]] \bigrootedtree[[[]]]
        & \bigrootedtree[] \bigrootedtree[] \bigrootedtree[[[]]]
        & \bigrootedtree[] \bigrootedtree[[]] \bigrootedtree[[]]
        & \bigrootedtree[] \bigrootedtree[] \bigrootedtree[] \bigrootedtree[[]]
        & \bigrootedtree[] \bigrootedtree[[],[[]]]
        & \bigrootedtree[] \bigrootedtree[[],[[]]]
        & \bigrootedtree[] \bigrootedtree[] \bigrootedtree[[],[]]
    \\[0.5cm]
    $p_\tau$
        & \bigrootedtree[]
        & \bigrootedtree[[]]
        & \bigrootedtree[[],[]]
        & \bigrootedtree[[],[]]
        & \bigrootedtree[[[]]]
        & \bigrootedtree[[],[[]]]
        & \bigrootedtree[[[]],[]]
        & \bigrootedtree[[[]],[]]
        & \bigrootedtree[[[]],[[]]]
        & \revB{\bigrootedtree[[]]}
        & \revB{\bigrootedtree[[[]]]}
        & \bigrootedtree[[],[]]
        & \bigrootedtree[[],[[]]]
        & \bigrootedtree[[]]
        & \bigrootedtree[[]]
        & \bigrootedtree[[],[]]
    \\
    \bottomrule
\end{tabular}
\end{table}

\begin{table}
    \centering
  \caption{All ordered subtrees of the tree with level sequence $[1,2,3,2,3]$.
           \revC{Each graph in the last row is an ordered subtree.
           The corresponding entry in the second row is the associated forest obtained
           when the subtree and adjacent edges are removed.  The top row shows the original
           tree, with the ordered subtree in black, removed edges as dashes, and the nodes
           and edges of the associated forest in red.}}
  \label{tbl:subtrees}
  \begin{tabular}{ccccccccccccccccc}
    \toprule
    $s$ & \bigrootedtree[[[]],[[]]]
        & \bigrootedtree[[[]],[[.,red node,edge=dashed]]]
        & \bigrootedtree[[[.,red node, edge=dashed]],[[]]]
        & \bigrootedtree[[[]],[.,red node,edge=dashed[.,red node,edge=red]]]
        & \bigrootedtree[[.,red node,edge=dashed[.,red node,edge=red]],[[]]]
        & \bigrootedtree[[[.,red node, edge=dashed]],[[.,red node, edge=dashed]]]
        & \bigrootedtree[[[.,red node,edge=dashed]],[.,red node,edge=dashed[.,red node,edge=red]]]
        & \bigrootedtree[[.,red node,edge=dashed[.,red node,edge=red]],[[.,red node,edge=dashed]]]
        & \bigrootedtree[[.,red node,edge=dashed[.,red node,edge=red]],[.,red node,edge=dashed[.,red node,edge=red]]]
        & \bigrootedtree[.,red node[.,red node,edge=red[.,red node,edge=red]],[.,red node,edge=red[.,red node,edge=red]]]
    \\[0.5cm]
    $\tau \setminus s$ & $\varnothing$
        & \bigrootedtree[]
        & \bigrootedtree[]
        & \bigrootedtree[[]]
        & \bigrootedtree[[]]
        & \bigrootedtree[] \bigrootedtree[]
        & \bigrootedtree[] \bigrootedtree[[]]
        & \bigrootedtree[] \bigrootedtree[[]]
        & \bigrootedtree[[]] \bigrootedtree[[]]
        & \bigrootedtree[[[]],[[]]]
    \\[0.5cm]
    $s_\tau$ & \bigrootedtree[[[]],[[]]]
        & \bigrootedtree[[[]],[]]
        & \bigrootedtree[[[]],[]]
        & \bigrootedtree[[[]]]
        & \bigrootedtree[[[]]]
        & \bigrootedtree[[],[]]
        & \bigrootedtree[[]]
        & \bigrootedtree[[]]
        & \bigrootedtree[]
        & $\varnothing$
    \\
    \bottomrule
\end{tabular}
\end{table}

\end{landscape}

\subsubsection{RK22 modified equation}

As an example, the modified equation for the method \eqref{rk22} applied to
a generic ODE \eqref{ode} is obtained in \bseries by
\begin{verbatim}
    using BSeries, Latexify, SymPy
    α = symbols("α", real=true)
    A = [0 0; 1/(2α) 0]; b = [1-α, α]; c = [0, 1/(2α)]
    coefficients = modified_equation(A, b, c, 4)
    println(latexify(coefficients, reduce_order_by=1, cdot=false))
\end{verbatim}
Note that we set \verb|reduce_order_by=1| to obtain the leading terms of the
B-series of $f_h$ instead of $h f_h$ as follows:
\begin{multline}
  f_h =
  F_{f}\mathopen{}\left( \rootedtree[] \right)\mathclose{}
  + \frac{ - h^{2}}{6} F_{f}\mathopen{}\left( \rootedtree[[[]]] \right)\mathclose{}
  + h^{2} \left( \frac{-1}{6} + \frac{1}{8 \alpha} \right) F_{f}\mathopen{}\left( \rootedtree[[][]] \right)\mathclose{}
  + \frac{h^{3}}{8} F_{f}\mathopen{}\left( \rootedtree[[[[]]]] \right)\mathclose{}
  + h^{3} \left( \frac{1}{8} - \frac{1}{16 \alpha} \right) F_{f}\mathopen{}\left( \rootedtree[[[][]]] \right)\mathclose{}
  \\
  + h^{3} \left( \frac{1}{8} - \frac{1}{8 \alpha} \right) F_{f}\mathopen{}\left( \rootedtree[[[]][]] \right)\mathclose{}
  + h^{3} \left( \frac{1}{24} - \frac{1}{16 \alpha} + \frac{1}{48 \alpha^{2}} \right) F_{f}\mathopen{}\left( \rootedtree[[][][]] \right)\mathclose{}
  + \order(h^4).
\end{multline}

\subsubsection{Example: Lotka-Volterra model}

As an example of greater mathematical interest,
we study the Lotka-Volterra model for population dynamics:
\begin{subequations} \label{LV}
\begin{align}
\dot{p}(t) & = (2-q)p, \\
\dot{q}(t) & = (p-1)q.
\end{align}
\end{subequations}
A backward error analysis of the explicit Euler method
\begin{align}
    y_{n+1} & = y_n + f(y)
\end{align}
applied to \eqref{LV} is presented in \cite[p. 340]{hairer2013geometric}, considering
only the first correction term (of order $h$).  Here we demonstrate the effectiveness
of {\bseries} by conducting an analysis including much higher order terms.

\begin{figure}
  \subfloat[$h=0.1, T=10$\label{fig:LVh1}]{
        \includegraphics[width=0.49\textwidth]{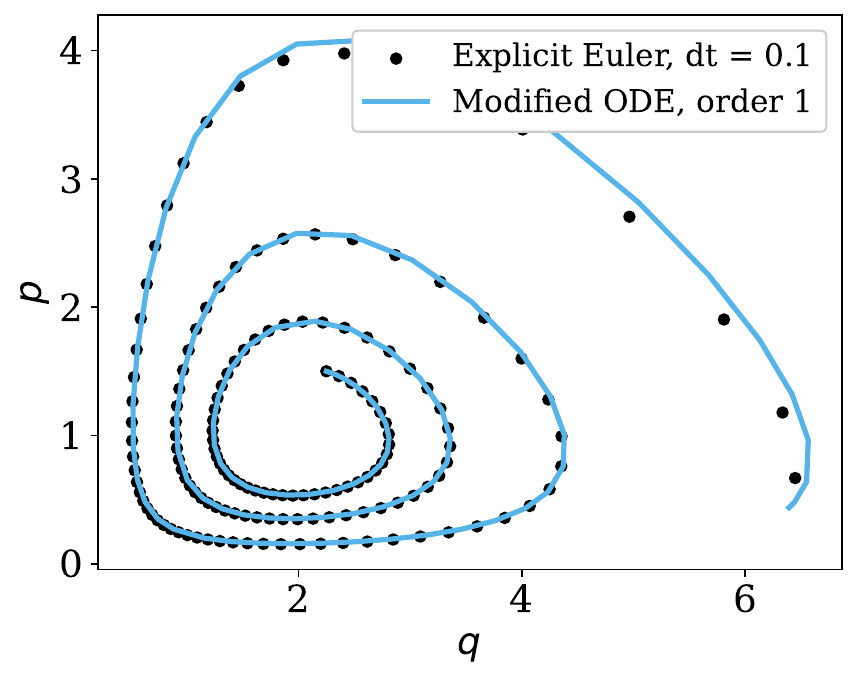}}
  \subfloat[$h=0.11, T=10$\label{fig:LVh2}]{
        \includegraphics[width=0.49\textwidth]{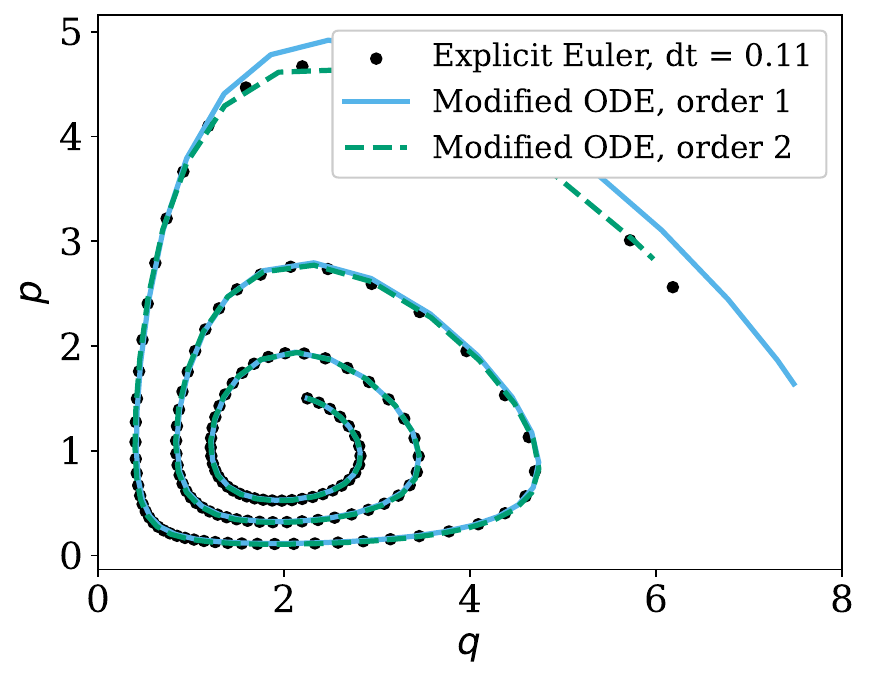}} \\
  \subfloat[$h=0.12, T=10$\label{fig:LVh3}]{
        \includegraphics[width=0.49\textwidth]{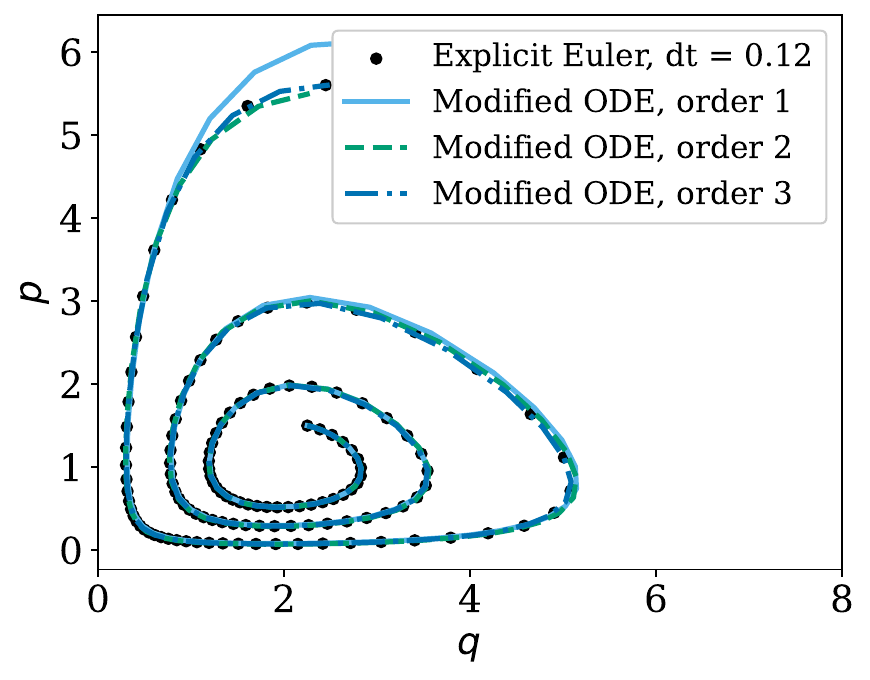}}
  \subfloat[$h=0.1$, $T=66.4$\label{fig:LVh4}]{
        \includegraphics[width=0.49\textwidth]{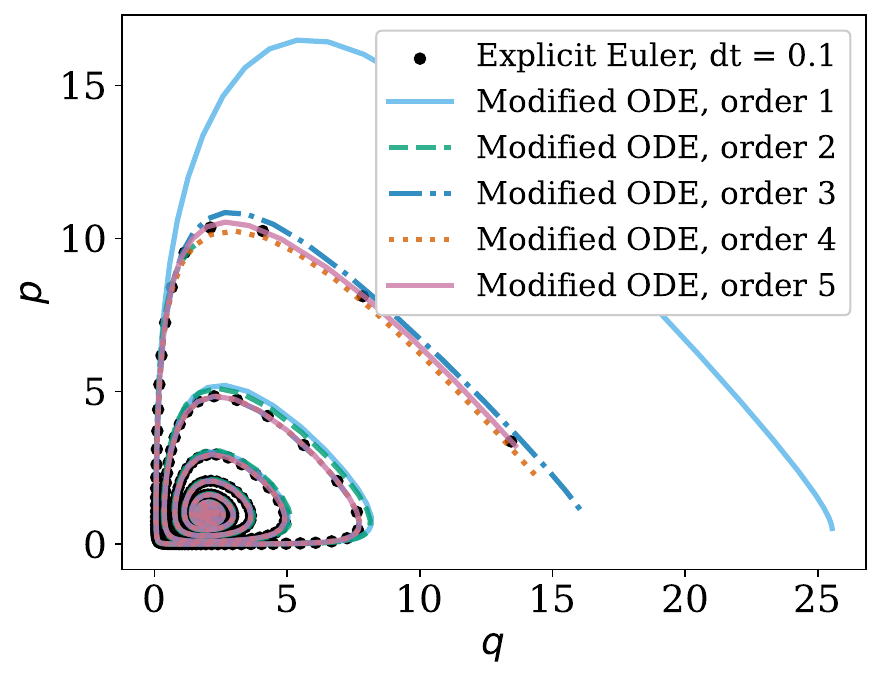}}
  \caption{Numerical solutions of the Lotka-Volterra system with final time $T$
           by the explicit Euler method with fixed time step size $h$ and
           reference solutions of corresponding modified equations.}
\end{figure}

We can set up the ODE system as follows:
\begin{verbatim}
    using BSeries, Latexify, SymPy
    h = symbols("h", real=true)
    u = p, q = symbols("p, q", real=true)
    f = [p * (2 - q), q * (p - 1)]
\end{verbatim}
After setting $A, b$ to the coefficients of the explicit Euler method, we can compute
the modified equation to any desired order with the command:
\begin{verbatim}
    A = fill(0//1, 1, 1); b = [1//1]; c = [0//1]
    series = modified_equation(f, u, h, A, b, c, 2)
\end{verbatim}
Computing just the $\order(h)$ terms yields
\begin{align} \label{LVh1}
\dot{p}(t) & \approx \frac{p \left(h \left(q \left(p - 1\right) - \left(q - 2\right)^{2}\right) - 2 q + 4\right)}{2}  \\
\dot{q}(t) & \approx \frac{q \left(h \left(p \left(q - 2\right) - \left(p - 1\right)^{2}\right) + 2 p - 2\right)}{2},
\end{align}
which agrees with what is presented in \cite{hairer2013geometric}.  Solving \eqref{LVh1}
to high accuracy yields a result that agrees well with the numerical solution from
the explicit Euler method when $h$ and $t$ are not too large; see Figure \ref{fig:LVh1}.
In Figure \ref{fig:LVh2} we show results obtained with a slightly larger step size; now the first-order modified equation
becomes inaccurate but the second-order equation still gives an accurate result on this scale.  If the step size
or final time is taken larger, more terms are required in order to obtain a modified equation that
agrees with the numerical solution, as shown in Figures \ref{fig:LVh3} and \ref{fig:LVh4}.
Even for small step sizes, several terms may be required to obtain good accuracy over longer times,
as shown in Figure \ref{fig:LVh4}, where we have included terms up to $\order(h^5)$.

%\subsubsection{Example: nonlinear pendulum}
%Reproduce examples from \url{https://nbviewer.jupyter.org/gist/ketch/28d83ec4134e62f8bd8ec4b3b6cccc4a}.

\subsection{Modifying integrators\label{sec:modifying-integrators}}
In the previous section we saw that we could accurately determine the behavior
of the explicit Euler method for the Lotka-Volterra system by deriving
increasingly accurate modified equations.  Nevertheless, the method's
qualitative behavior remains completely different from that of the true
solution (which, for the prescribed conditions, is periodic).  In this
section, we instead seek to modify the right-hand-side $f(y)$ to obtain a
perturbation $f_h(y)$ so that
when the explicit Euler method is applied to $f_h(y)$ it
yields a solution that matches the exact solution of the original problem $f(y)$.
This approach is referred to as a \emph{modifying integrator}; i.e. the integrator
modifies the ODE.

In terms of B-series, we want to obtain the exact solution B-series when we
substitute the modified RHS into our numerical method:
$$
    B(\coeff, B(\coeffb,hf,\cdot),y) = B(e, hf, y)
$$
Here $B(\coeff,hf,y)$ is the B-series of our numerical method (in this case, explicit Euler)
while $B(\coeffb,hf,y)$ is the B-series of the modified equation.  We see that we have
$(\coeffb \star \coeff)(\tau) = e(\tau)$.

\begin{figure}
  \subfloat[Exact solution and explicit Euler solution.\label{fig:MI1}]{
        \includegraphics[width=0.49\textwidth]{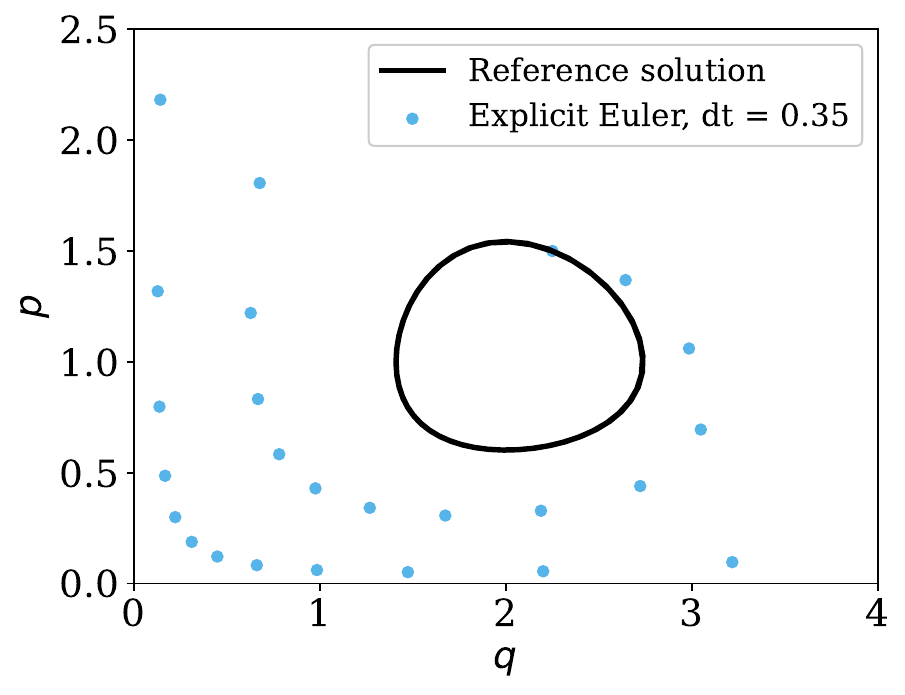}}
  \subfloat[Solutions of modifying integrators of orders 1-3.\label{fig:MI2}]{
        \includegraphics[width=0.49\textwidth]{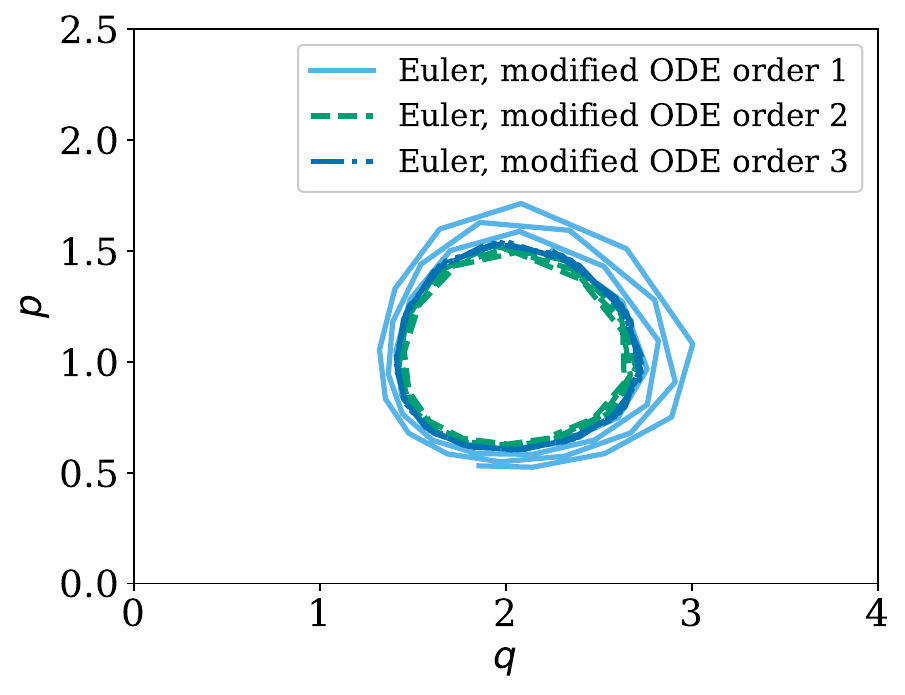}}
  \caption{Numerical solutions of the modifying integrator system
            corresponding to the Lotka-Volterra equations and explicit
            Euler method.\label{fig:modint}}
\end{figure}

Results for the case at hand are shown in Figure \ref{fig:modint}.
The left figure shows the exact solution and the numerical solution
of the original problem, while the right figure shows solutions of
the modifying integrator systems of orders 1, 2, and 3.  We see that
the first-order modifying integrator system still exhibits growth,
though it is greatly reduced compared to the Euler solution of the original
system.  The second-order modifying integrator is actually dissipative.
The third-order modifying integrator is indistinguishable from the exact solution
at this scale and is nearly periodic.

\subsection{A special RK method and ODE pair} \label{sec:nlo}
In this section we study the application of the explicit midpoint Runge-Kutta method
\begin{subequations} \label{midpoint}
\begin{align}
    Y_2 & = y_n + \frac{h}{2}f(y_n), \\
    y_{n+1} & = y_n + h f(Y_2),
\end{align}
\end{subequations}
applied to the nonlinear oscillator problem
\begin{subequations} \label{nlo}
\begin{align}
\dot{p}(t) & = \frac{-q}{p^2 + q^2}, \\
\dot{q}(t) & = \frac{ p}{p^2 + q^2}.
\end{align}
\end{subequations}
Solutions of \eqref{nlo} are periodic, forming circles
in the $p-q$ plane, but with period depending on the initial
energy \revB{$\beta^2 = \|y(0)\|^2$}, which is conserved:
\begin{align} \label{num_sol}
    y(t) & = \beta \begin{bmatrix} \cos(t/\beta^2) \\ \sin(t/\beta^2) \end{bmatrix}.
\end{align}

\revB{The numerical solution obtained by application of \eqref{midpoint} to
\eqref{nlo} turns out to be
\begin{align} \label{numsol}
y_n = \beta \begin{bmatrix} \cos(n\theta(h/\beta^2)) \\ \sin(n\theta(h/\beta^2)) \end{bmatrix},
\quad \text{where }
    \theta(\mu) = \arccos\left( \frac{1-\frac{\mu^2}{4}}{1+\frac{\mu^2}{4}} \right).
\end{align}}
This can be verified by substituting \eqref{numsol}
into \eqref{midpoint} and solving for $\theta$, noting that
the given value solves both of the resulting equations.
As studied in \cite[Section~5.1]{2019_energystability}, the midpoint
method \eqref{midpoint} has the remarkable property that it conserves
the energy exactly for \eqref{nlo} \emph{for any step size}.  As far
as we are aware, this behavior is unique among explicit discretizations
of initial value ODEs.  The backward error analysis of this example was
our initial motivation for developing \bseries.

\subsubsection{Modified equation analysis}
The modified equation (up to order eight) for this example can be obtained with the following code:
\begin{verbatim}
    using BSeries, StaticArrays, SymPy

    A = @SArray [0 0; 1//2 0]
    b = @SArray [0, 1//1]
    c = @SArray [0, 1//2]

    u = symbols("p,q")
    h = symbols("h")
    f = [-u[2], u[1]] / (u[1]^2 + u[2]^2)

    modified_equation(f, u, h, A, b, c, 8);
\end{verbatim}
After some cleaning up, we find that the modified equation RHS is closely
related to the original RHS: \revB{
\begin{align} \label{oh6modeq}
    f_h(y) = f(y) g(h/\|y\|^2)
\end{align}
\revB{where $g$ is a function of one variable, with the first terms in its power series (obtained from
the code above) being
\begin{align} \label{pofz}
    g(z) = 1 - \frac{z^2}{12} + \frac{z^4}{80} - \frac{z^6}{448} + \frac{z^8}{2304} + \order(z^{10}).
\end{align}}}
In fact, this structure could be anticipated because
the Hamiltonian structure of \eqref{nlo} implies that the even-order derivatives of
$f$ are proportional to $f$, while odd-order derivatives are proportional to $y$.
This means that any odd-order derivative terms would be non-conservative, whereas
(since we know the midpoint method is energy-conservative)
\cite[Corollary~IX.5.3]{hairer2013geometric} guarantees that every truncation of the
modified equation B-series must be conservative.

\revB{For this example we can find $f_h$ in an alternative and more direct
way\footnote{This was pointed out in a private communication
by Prof. Ernst Hairer.}.
From \eqref{num_sol} we have
\begin{align*}
    y_h(h) & =  \beta \begin{bmatrix} \cos(\theta(h/\beta^2)) \\ \sin(\theta(h/\beta^2)) \end{bmatrix}
\end{align*}
and differentiating both sides with respect to $h$ yields
\begin{align*}
y_h'(h) & =  \beta  \begin{bmatrix} -\sin(\theta(h/\beta^2)) \\ \cos(\theta(h/\beta^2)) \end{bmatrix} \frac{d}{dh}\theta(h/\beta^2)\\
        & = f(y_h(h)) \theta'(h/\beta^2).
\end{align*}
Comparing with \eqref{oh6modeq}, we see that
\begin{align} \label{dirmodeq}
    g(z) & = \frac{d}{dz} \arccos\left( \frac{1-\frac{z^2}{4}}{1+\frac{z^2}{4}} \right) = 2 \arctan(z/2).
\end{align}}

This also implies that the midpoint method is symmetric (in the sense of \cite{hairer2013geometric})
for this problem.  This can be verified directly by checking that \eqref{num_sol}
is the solution obtained with the adjoint of the midpoint method, which is
\begin{subequations} \label{rk2adj}
\begin{align}
    Y_1 & = y_n + h f(Y_1) - \frac{h}{2} f(Y_2) \\
    y_{n+1} = Y_2 & = y_n + h f(Y_1).
\end{align}
\end{subequations}

It is important to note that the explicit midpoint method is \emph{not}
symmetric for general ODEs; the B-series for \eqref{midpoint} for a
general RHS $f$ includes both even- and odd-order derivatives.  Thus
the behavior observed in this example results from a (seemingly very
fortuitous) cancellation of terms at every odd order in $h$.  At order
$h^5$, there are 16 (out of 20) coefficients in the modified equation
for \eqref{midpoint} that are non-zero, but their sum when multiplied
by the elementary differentials corresponding to \eqref{nlo} is zero.
The method and the problem appear to possess a complementary structure.

\subsubsection{Modifying integrator}
We can also construct a modifying integrator for this problem.
Again using \bseries we find that, up to terms of order $h^4$, we have
$$
    \hat{f}_h(y) = f(y) \hat{g}_4(h/\|y\|^2),
$$
where $f(y)$ is the right hand side of the original
problem \eqref{nlo} and
$$
    \hat{g}_4(z) := 1 + \frac{z^2}{12} + \frac{z^4}{20}.
$$
We see that the right hand side is again Hamiltonian at each order, and we might conjecture that $\hat{f}_h$ will have this structure at every order.  However, continuing the calculation up to terms of order $h^7$ yields
\begin{equation}
\label{mod_int_nlo_midpoint}
    \hat{f}_h(y) = f(y)\hat{g}_6(h/\|y\|^2) + y
    \left(\frac{h^5}{48\|y\|^{12}} +\frac{31 h^7}{640 \|y\|^{16}}\right)
\end{equation}
where
$$
    \hat{g}_6(z) := \hat{g}_4(z) + \frac{127}{2016} z^6.
$$
We see that the truncations of the modifying integrator
system are not Hamiltonian (for order $h^5$ and higher).

\revC{%
The energy and error of numerical solutions of the explicit
midpoint method \eqref{midpoint} applied to the nonlinear
oscillator \eqref{nlo} and the corresponding modifying integrator
equation \eqref{mod_int_nlo_midpoint} for two choices of the
constant time step size $h$ are shown in
Figure~\ref{fig:nonlinear_osc_mod_int}.
Initially, the modifying integrator approach yields a much more
accurate solution (since it results in an eighth-order approximation
instead of a second-order one).
However, the error growth rate in time is larger since the energy is
not conserved. Thus, the energy-conserving explicit midpoint method
can have a smaller error for large times. However, this asymptotic
regime depends on the chosen time step size. If $h$ is small enough,
the increased order of accuracy of the modifying integrator approach
makes it advantageous until both methods yield approximately 100\%
error.

\begin{figure}
  \subfloat[Error, $h=0.6$]{
        \includegraphics[width=0.49\textwidth]{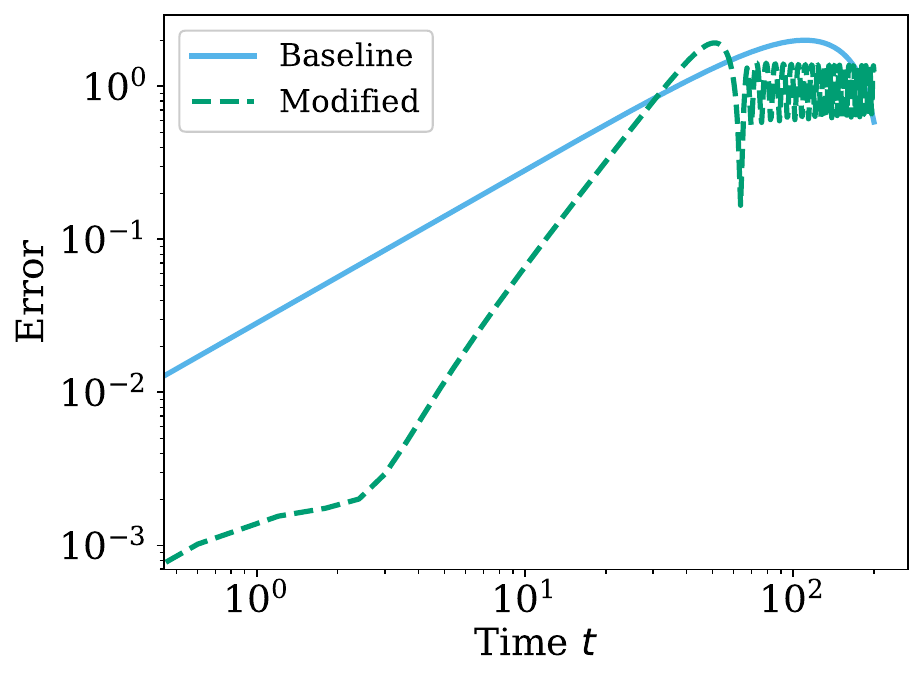}}
  \subfloat[Energy, $h=0.6$]{
        \includegraphics[width=0.49\textwidth]{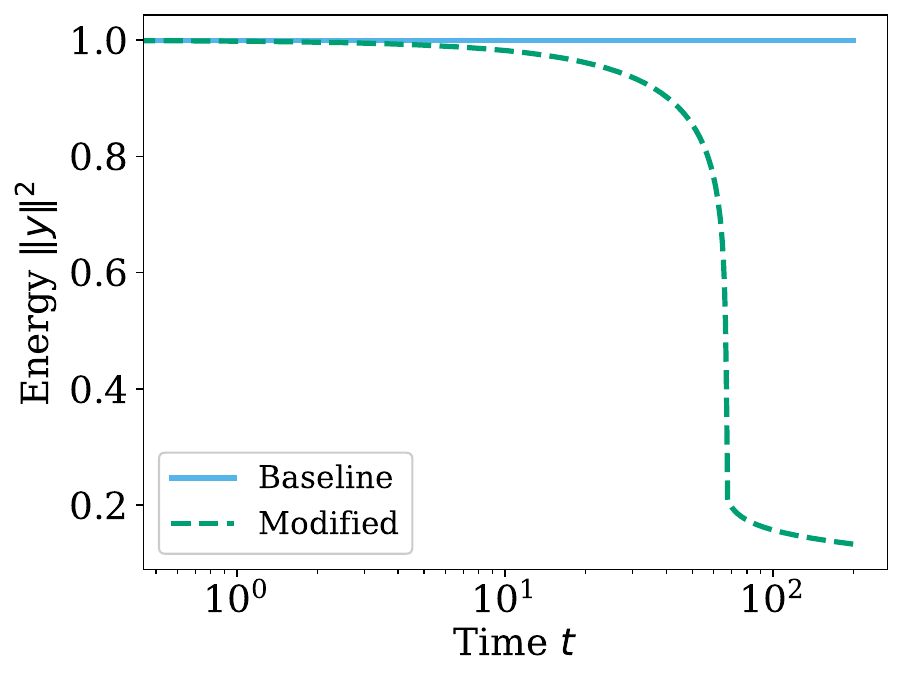}} \\
  \subfloat[Error, $h=0.3$]{
        \includegraphics[width=0.49\textwidth]{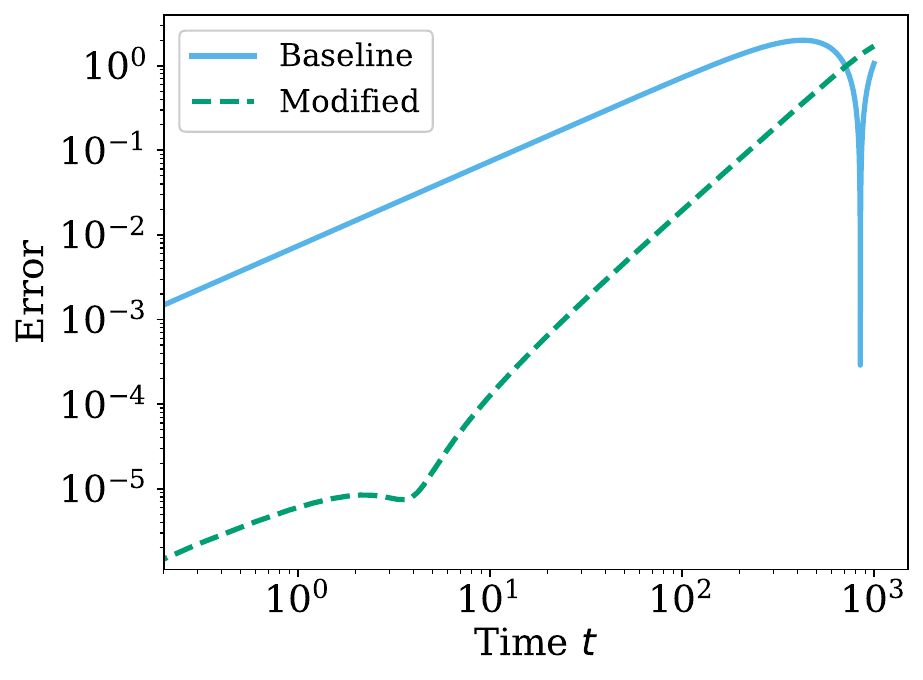}}
  \subfloat[Energy, $h=0.3$]{
        \includegraphics[width=0.49\textwidth]{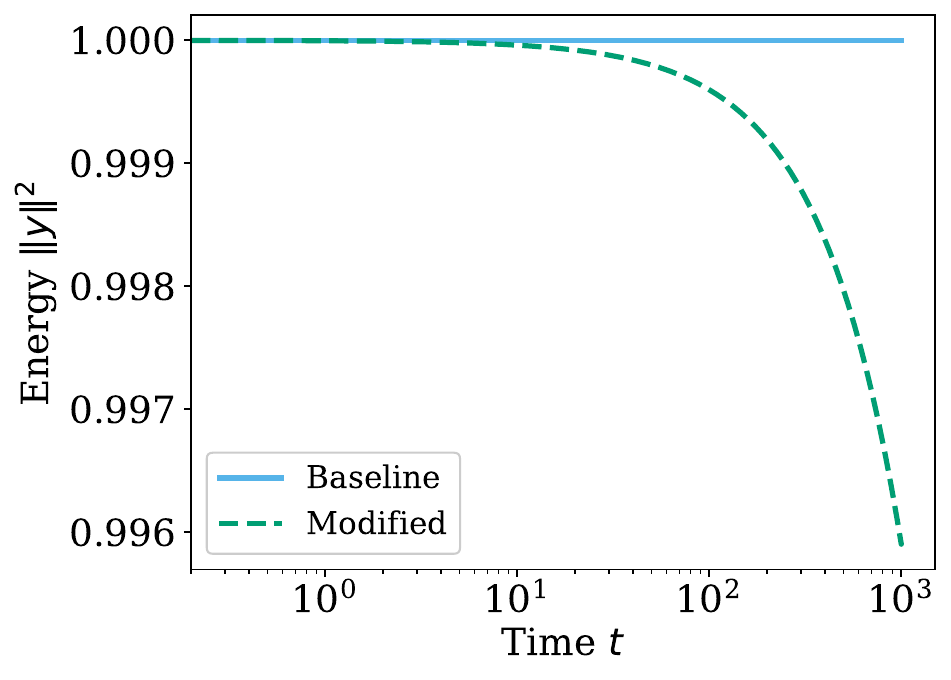}}
  \caption{\revC{Error and energy of numerical solutions of the
           explicit midpoint method \eqref{midpoint} with
           fixed time step size $h$ for
           the nonlinear oscillator \eqref{nlo} and the corresponding
           modifying integrator equation (including terms of
           order $h^7$).}}
  \label{fig:nonlinear_osc_mod_int}
\end{figure}
} % \revC

\revA{
\section{B-series for other classes of methods}
The historical development of B-series is closely intertwined with that
of Runge-Kutta methods, and Runge-Kutta methods (and their generalizations)
continue to be the primary area of application for B-series.  Nevertheless,
B-series are also used in the analysis of other classes of methods, and
\bseries can be used for this purpose.  In Section \ref{sec:ARK} we discuss
application to additive Runge-Kutta methods.  In Section \ref{sec:AVF}
we present the construction and
analysis of the B-series and modified equation of the average vector field
(AVF) method \cite{mclachlan1999geometric}.

\revA{%
\subsection{Additive Runge-Kutta methods}\label{sec:ARK}

Many time integration methods outside the class of Runge-Kutta methods
-- such as partitioned RK methods \cite[Section~III.2]{hairer2013geometric} and
Runge-Kutta-Nyström methods \cite[Section~II.14]{hairer1993} --
can be represented using generalizations of B-series based on colored
rooted trees.
Colored rooted trees are also called $N$-trees \cite{araujo1997symplectic}
or $P$-trees \cite{hairer1981order}. The corresponding generalized
B-series are also called NB-series \cite{araujo1997symplectic} or P-series
\cite{hairer1981order}, depending on the context.
Here, we briefly discuss the extension of B-series to additive
Runge-Kutta (ARK) methods (a superset of partitioned RK methods) and
functions on trees to colored rooted trees, following
\cite{araujo1997symplectic}.

For an additive ODE
\begin{equation} \label{additive_ode}
  \dot y(t) = \sum_{\nu=1}^N f^\nu\bigl( y(t) \bigr)
\end{equation}
with $N \in \mathbb{N}$, an additive Runge-Kutta method performs a step from
$y_{n}$ to $y_{n+1}$ via
\begin{subequations} \label{ARK}
\begin{align}
    Y_i & = y_n + h \sum_{\nu=1}^N \sum_{j=1}^s a_{ij}^{\nu} f^{\nu}(Y_j), \\
    y_{n+1} & = y_n + h \sum_{\nu=1}^N \sum_{j=1}^s b_j^{\nu} f^{\nu}(Y_j).
\end{align}
\end{subequations}
The Butcher coefficients are $A^{\nu} = (a_{ij}^{\nu})_{i,j=1}^{s}$,
$b^{\nu} = (b_j^{\nu})_{j=1}^s$, $\nu \in \{1, \dots, N\}$. Here, we assumed
again an autonomous form of the ODE. The same analysis applies to ARK methods
satisfying the row-sum condition.

We study an example of an additive ODE and ARK method in Section 
\ref{sec:sv}.  But first we briefly discuss the extension of the theory
of rooted trees and B-series to such objects.

\subsubsection{Colored rooted trees}
The solution map corresponding to an ARK method can be represented by
an infinite sum whose terms are in one-to-one correspondence with colored
rooted trees.
Each node in a colored rooted tree has an associated color, which
can also be represented by an integer $\nu \in \{1, 2, \dots, N\}$
corresponding to a term $f^\nu$ in \eqref{additive_ode}.
Each tree has an associated elementary differential, which is identical
to that of the corresponding uncolored tree except that each factor involving
a derivative of $f$ now involves the corresponding $f^\nu$.
For example, let $N=2$ and (for clarity) write $f^1 = f$, $f^2=g$.
Some examples of elementary differentials are as follows:}
\begin{align*}
F\left(\bigrootedtree[[[.,gray node]]]\right) & = f^j_k f^k_l g^l &
F\left(\bigrootedtree[.,gray node [[] [.,gray node] ]]\right) & = g^j_k f^k_{lm} f^l g^m.
\end{align*}
\revA{
In a similar way, the elementary weight of a colored rooted tree is
identical to that of the uncolored tree except that the coefficients
$A, b$ associated with each node are replaced by $A^\nu, b^\nu$ where $\nu$
is the color of the node.  Some examples of elementary weights for
bicolored trees are as follows:}
\begin{align*}
\Phi\left(\bigrootedtree[[[.,gray node]]]\right) & = \sum_{j,k,l} b^1_j a^1_{jk} a^2_{kl} &
\Phi\left(\bigrootedtree[.,gray node [[] [.,gray node] ]]\right) & = \sum_{j,k,l,m} b^2_j a^1_{jk} a^1_{kl} a^2_{km}.
\end{align*}
\revA{
In the code, we represent colored rooted trees using a color sequence
in addition to the level sequence. The algorithms such as the generation
of the partitions are extended correspondingly. The generation of all colored
rooted trees extends the algorithm of \cite{beyer1980constant} by
iterating over all level sequences at the outer level and all color sequences
in an inner iteration; this inner iteration needs to ensure that the color
sequence yields a canonical representation, which is extended to colored rooted
trees by additional lexicographic comparison of the color sequences after
the usual lexicographic comparison of the level sequences.

\subsubsection{Störmer-Verlet method as additive Runge-Kutta method} \label{sec:sv}

Consider the Kepler problem
\begin{equation}
  \dot q(t) = p,
  \quad
  \dot p(t) = g(q) = \frac{q}{|q|^3},
\end{equation}
with momentum $p = (p_1, p_2)$ and position $q = (q_1, q_2)$. This is an example
of an additive ODE \eqref{additive_ode} with $N = 2$ and
\begin{equation}
  f^1(q, p) = \begin{pmatrix} p \\ 0 \end{pmatrix},
  \quad
  f^2(q, p) = \begin{pmatrix} 0 \\ q / |q|^3 \end{pmatrix}.
\end{equation}
The classical Störmer-Verlet method
\begin{equation}
\begin{aligned}
  p_{n+1/2} &= p_{n} + \frac{h}{2} g(q_{n}), \\
  q_{n + 1} &= q_{n} + h p_{n+1/2}, \\
  p_{n + 1} &= p_{n+1/2} + \frac{h}{2} g(q_{n+1}),
\end{aligned}
\end{equation}
can be interpreted as an ARK method \eqref{ARK} with $N = 2$ and coefficients
\begin{equation}
  A^1 = \begin{pmatrix} 0 & 0 \\ 1/2 & 1/2 \end{pmatrix},
  b^1 = \begin{pmatrix} 1/2 \\ 1/2 \end{pmatrix},
  \quad
  A^2 = \begin{pmatrix} 1/2 & 0 \\ 1/2 & 0 \end{pmatrix},
  b^2 = \begin{pmatrix} 1/2 \\ 1/2 \end{pmatrix},
\end{equation}
cf. Table~II.2.1 of \cite{hairer2013geometric}. Although these coefficients
look like those of an implicit method, they result in the explicit Störmer-Verlet
method since $f^1(q, p)$ influences only $q$ and does not depend on $q$ while
$f^2(q, p)$ influences only $p$ and does not depend on $p$.

Computing the modified equation up to terms of order $h^7$ using \bseries
reveals that there are no terms proportional to odd powers of $h$, in accordance
with the well-known symplectic property of the method.

Next we consider improving the accuracy of the Störmer-Verlet method for Kepler's problem
by applying the modifying integrator approach, using \bseries, up to
order four.  The resulting modifying integrator equations include the following equation for $q_1$:
\begin{equation}
  \dot q_1(t) = p_{1} +
  \frac{h^{2}}{6} \left( p_{1} \left( \frac{3 q_{1}^{2}}{\left( q_{1}^{2} + q_{2}^{2} \right)^{\frac{5}{2}}} - \frac{1}{\left( q_{1}^{2} + q_{2}^{2} \right)^{\frac{3}{2}}} \right) + \frac{3 p_{2} q_{1} q_{2}}{\left( q_{1}^{2} + q_{2}^{2} \right)^{\frac{5}{2}}} \right).
\end{equation}
We see that there are no third-order terms (as expected for this symmetric
method). However, since $q$ appears on the right hand side, the method is no
longer explicit.  We conclude that this does not appear to be an appealing
approach.  The code for this example is accessible in our reproducibility
repository \cite{ketcheson2021computingRepro}.
} % \revA

\subsection{The average vector field method}\label{sec:AVF}

The AVF method takes the form
\begin{align}
    y_{n+1} = y_{n} + h \int_0^1 f\bigl(\xi y_{n+1} + (1 - \xi) y_{n}\bigr) \mathrm{d} \xi.
\end{align}
It is a B-series method, with coefficients given explicitly as
\begin{align}
    b(\rootedtree[]) &= 1, \\
    b([t_1, ..., t_n]) &= b(t_1)...b(t_n) / (n + 1);
\end{align}
see \cite{quispel2008new,celledoni2009energy}.
We can instantiate this method in \bseries with the following code:}
\begin{verbatim}
using BSeries

series = bseries(5) do t, series
    if order(t) in (0, 1)
        return 1 // 1
    else
        v = 1 // 1
        n = 0
        for subtree in SubtreeIterator(t)
            v *= series[subtree]
            n += 1
        end
        return v / (n + 1)
    end
end
\end{verbatim}
\revA{Here we have generated the B-series up to order five.
The energy-preserving property of the method can be checked by considering
its modified equation \cite{celledoni2010energy}, which we do as follows:}
\begin{verbatim}
meq = modified_equation(series)
println(latexify(meq, reduce_order_by=1, cdot=false))
\end{verbatim}
\revA{which yields the modified equation with}
\begin{align*}
f_h = F_{f}\mathopen{}\left( \rootedtree[] \right)\mathclose{} + \frac{1}{12} h^{2} F_{f}\mathopen{}\left( \rootedtree[[[]]] \right)\mathclose{} + \frac{1}{80} h^{4} F_{f}\mathopen{}\left( \rootedtree[[[[[]]]]] \right)\mathclose{} + \frac{1}{180} h^{4} F_{f}\mathopen{}\left( \rootedtree[[[[][]]]] \right)\mathclose{} + \frac{1}{360} h^{4} F_{f}\mathopen{}\left( \rootedtree[[[[]][]]] \right)\mathclose{} + \frac{-1}{180} h^{4} F_{f}\mathopen{}\left( \rootedtree[[[[]]][]] \right)\mathclose{} \\
+ \frac{-1}{720} h^{4} F_{f}\mathopen{}\left( \rootedtree[[[][][]]] \right)\mathclose{} + \frac{-1}{180} h^{4} F_{f}\mathopen{}\left( \rootedtree[[[][]][]] \right)\mathclose{} + \frac{1}{360} h^{4} F_{f}\mathopen{}\left( \rootedtree[[[]][[]]] \right)\mathclose{} + \frac{-1}{720} h^{4} F_{f}\mathopen{}\left( \rootedtree[[[]][][]] \right)\mathclose{}.
\end{align*}
\revA{This reproduces the formula given at the top of p. 679 of \cite{celledoni2010energy},
demonstrating the energy-preserving property (at least up to the order $h^4$ terms).
It is straightforward to study higher-order terms in the same way.}

\section{Implementation and performance}

We have developed \bseries and RootedTrees.jl together to get reasonably efficient
libraries that allow working with high-order representations of B-series. Of course,
the increasing number of elementary differentials of a given order
makes explicit calculations infeasible beyond some order.  We tried to
make this order as high as reasonably possible without increasing the code complexity
too much.  In this section we explain some of the design choices and optimizations
that have been made.  Additional performance optimizations that could be
\revC{implemented} in the future include parallel computation and memoization, which
which could increase the runtime efficiency of some parts at the cost of a higher
memory footprint and/or higher code complexity.

Before discussing some of the design decisions of \bseries related to computational
performance, we compare its efficiency to \pybs \cite{sundklakk2015pybs}.
Since \pybs uses lazy representations of B-series up to an arbitrary
order, we need to instantiate the coefficients for a reasonable comparison. Here,
we set up the B-series of the explicit midpoint method with Butcher tableau
\begin{equation}
  \begin{array}{c|cc}
  0 & 0 & \\
  \frac{1}{2} & \frac{1}{2} & 0 \\[0.1cm]
  \hline
  & 0 & 1
  \end{array}
\end{equation}
and sum the coefficients of the B-series of its modified equation up to order nine.
On an Intel® Core™ i7-8700K from 2017, this takes approximately
\SI{53}{\ms} with \bseries and
\SI{3.7}{\s} with \pybs.
Thus, \bseries is approximately 70x faster than \pybs for this calculation.
All source code and instructions required to reproduce these results are
available in our reproducibility repository \cite{ketcheson2021computingRepro}.

The first design decision of \bseries is the choice of programming language.
We chose Julia \cite{bezanson2017julia} since it provides convenient high-level
abstractions, interactive workflows, and efficient computations enabled via
compilation to native machine code. In fact, Julia can be at least as fast as
low-level code written in traditional programming languages such as
Fortran, C, and C++ \cite{elrod2021roadmap,ranocha2022adaptive}.

\subsection{Representation of rooted trees}

Many performance-related design decisions of \bseries are based on the
representation of rooted trees. Since RootedTrees.jl uses canonical level sequences
stored in arrays for this purpose, we adopted the same convention when extending
it with the functionality required for \bseries.

This convention is convenient when iterating over all rooted trees of a given
order since we can use the algorithm of \cite{beyer1980constant} directly.
Moreover, this representation is memory efficient with good data locality.
Additionally, it simplifies the computation of certain functions on rooted trees
such as its symmetry $\sigma$ since identical subtrees are always adjacent to
each other.

On the other hand, implementing the composition $\cdot$ or the substitution $\star$
using level sequences requires frequent memory movement. Moreover, the assumption
of a canonical representation requires that we sort level sequences of newly
created rooted trees, e.g., trees of the partition forest or the partition skeleton.

An alternative representation of rooted trees is based on linked lists of their
vertices. The advantages and disadvantages of such a representation and array-based
level sequences are basically mirrors of each other. We did not perform extensive
benchmarks of both possibilities. Instead, we implemented workarounds for most
common performance pitfalls of array-based level sequences.
For example, we
keep track of some internal buffers in RootedTrees.jl that allow us to avoid
dynamic memory allocations when sorting level sequences to obtain canonical
representations.
Some other performance optimizations are described below.

\subsection{Composition and substitution laws}

In the literature there exist different algorithmic approaches to the composition
and substitution laws, all in terms of rooted trees.  In our implementation we
have followed the work of Chartier et al. \cite{Chartier2010}.
Computing the coefficients of the composition of two B-series requires iterating
over all splittings of each tree. Similarly, the substitution law requires
iterating over all partitions of a rooted tree. Individually, each of these
operations requires iterating over the splitting forest or partition forest of
a tree.
To reduce the pressure on memory and keep data locality, we use lazy iterators
for these tasks. Similarly, we also use lazy iterators when iterating over all
rooted trees of a given order.

\revA{%
High-level pseudocode for finding all partitions of a tree (the key to the substitution law)
is given in Algorithms \ref{alg:partitions}--\ref{alg:pskeleton}, in order to illustrate
some of the performance considerations. These algorithms use lists storing the
partition forests/skeletons for clarity and simplicity or presentation; the efficient
implementations used for the substitution law use lazy iterators instead.
We iterate over binary sequences where each 1 corresponds to an edge that is kept
(solid edges in Table \ref{tbl:partitions}) and each 0 corresponds to an edge that
is removed (dashed edges in Table \ref{tbl:partitions}).
When executing the while-loops in Algorithms \ref{alg:pforest} and \ref{alg:pskeleton},
in which subtrees are removed from the initial tree, we work backwards from the end of
the level sequence.  Compared to beginning at the root, this reduces the amount of
memory moves on average since some nodes have already been removed.
Additionally, we use views into existing memory instead of copying parts of level
sequences to new arrays whenever possible.

\begin{algorithm}
\caption{Find all partition forests and corresponding skeletons of a tree $\tau$}\label{alg:partitions}
\begin{algorithmic}
    \State $m \gets |\tau|$
    \State forests $\gets$ []  \Comment{Empty list}
    \State partitions $\gets$ []
    \For{ edge\_set in (all binary sequences of length $m$)} \Comment{Edges marked zero are removed}
        \State Append partition\_forest($\tau$, edge\_set) to forests
        \State Append skeleton($\tau$, edge\_set) to skeletons
    \EndFor
\end{algorithmic}
\end{algorithm}

\begin{algorithm}
\begin{algorithmic}
\caption{Find the partition forest corresponding to a tree $\tau$ and subset of its edges}\label{alg:pforest}
\Procedure{partition\_forest}{$\tau$, edge\_set}
    \State forest $\gets$ [] \Comment{Empty list}
    \While{0 in edge\_set}
        \State $i \gets$ index of next edge marked 0
        \State remove from $\tau$ edge $i$ and its descendants
        \State append the removed subtree to forest
        \State delete entry $i$ from edge\_set
    \EndWhile
    \State append $\tau$ to forest
    \State \Return forest
\EndProcedure
\end{algorithmic}
\end{algorithm}

\begin{algorithm}
\begin{algorithmic}
\caption{Find the skeleton corresponding to a tree $\tau$ and subset of its edges}\label{alg:pskeleton}
\Procedure{skeleton}{$\tau$, edge\_set}
    \While{1 in edge\_set}
        \State $i \gets$ index of next edge marked 1
        \State Contract edge $i$ of $\tau$ by removing the child node and promoting its descendants
        \State delete entry $i$ from edge\_set
    \EndWhile
    \State \Return $\tau$
\EndProcedure
\end{algorithmic}
\end{algorithm}
} %\revA

The computation of modifying integrators benefits from additional checks
in the composition law: If the coefficient of the partition skeleton is zero,
there is no need to iterate over the associated partition forest, see
\eqref{substition}. This is particularly useful when computing the modified
integrator since we need to substitute the modifying integrator series into
the B-series of the numerical integrator.
By construction, many coefficients
of the B-series of explicit Runge-Kutta methods are zero since they involve
powers of the strictly lower diagonal matrix $A$. For example, computing the
coefficients of the modifying integrator of the explicit midpoint method up to
order nine takes approximately \SI{23}{\ms} while the same task for the modified
equation takes approximately \SI{52}{\ms}.

\subsection{Symbolic calculations}

Multiple symbolic computation packages are available in Julia, each with
certain strengths and weaknesses.
Therefore, we use Julia's efficient dispatch mechanism to support different
symbolic packages (SymPy.jl, Symbolics.jl, and SymEngine.jl) as backends in
\bseries. This allows users to choose a package based on their needs.
At the time of writing, we often use SymPy.jl for calculations that should be
rendered directly in \LaTeX, Symbolics.jl when we want to create efficient
Julia functions from symbolic computations, and SymEngine.jl when performance
matters most and we are mostly interested in symbolic expressions after inserting
specific initial values.

To compare the relative performance of different symbolic backends, we measured
the time to compute the B-series of the modifying integrator of the explicit
midpoint method \eqref{midpoint} specialized to the nonlinear oscillator
\eqref{nlo} up to order eight and to substitute the initial condition $(0, 1)$
into the resulting symbolic expressions. Using current versions of all packages
at the time of writing, SymEngine.jl is several orders of magnitude faster than
the other two alternatives; see Table~\ref{tab:performance-symbolics}.

\begin{table}[ht]
\centering
\caption{Micro-benchmarks of different symbolic backends supported by \bseries
          for computing the B-series of the explicit midpoint method specialized
          to the nonlinear oscillator discussed in Section~\ref{sec:nlo}.
          Note the different units of time (seconds and milliseconds).}
\label{tab:performance-symbolics}
\begin{tabular}{l S[table-format=1.6] S}
  \toprule
  \\
    & \multicolumn{1}{c}{Computing the B-series}
    & \multicolumn{1}{c}{Substituting the initial condition}
  \\
  \midrule
  SymPy.jl     & \SI{1.811 \pm 20.118e-3}{\s}    & \SI{212.317 \pm 9.885}{\ms}
  \\
  Symbolics.jl & \SI{710.232 \pm 21.727}{\ms} & \SI{2.565 \pm 2.848e-3}{\s}
  \\
  SymEngine.jl & \SI{15.768 \pm 7.733}{\ms}   & \SI{3.353 \pm 0.271}{\ms}
  \\
  \bottomrule
\end{tabular}
\end{table}

The performance optimizations for symbolic computations in \bseries are orthogonal
to the choice of specific backend. They are mostly related to programming approaches
and the use of information available from analytical theorems. For example, we
use a dynamic programming approach to determine elementary differentials of
symbolic functions, i.e., we compute all elementary differentials up to a specified
order upfront, allowing us to re-use lower order differentials when computing
higher order ones. In addition, we use Schwarz's theorem to reduce the number
of symbolic differentiation steps by making use of symmetry properties of
higher-order derivatives.

\section{Conclusion}
B-series are fascinating objects that lie at the intersection of continuous and discrete,
as well as pure and applied, mathematics.
The package \bseries makes it straightforward to apply B-series analysis to any Runge-Kutta method
and/or initial value ODE.  Most applications of B-series in the literature are limited to
lower order in $h$ due to the combinatorial explosion of terms that must be dealt with; \bseries
makes it possible to work effectively and efficiently with series up to much higher order in $h$,
thus automating backward error analysis for discretizations of nonlinear ODE
systems.  We hope that \bseries will also enable continued advances in our
understanding of numerical integrators.
Furthermore, since renormalization in quantum field theory is known also to rely on
the Hopf algebra structure implemented in \bseries
\cite{connes2000renormalization}, this package may be useful for calculations in that domain.

There are many potential areas for future extension of \bseries, corresponding to the
wide range of applications of B-series; these include B-series for partitioned methods,
multistep methods, multiderivative methods, exponential methods, and Rosenbrock methods.

\section*{Acknowledgements}

We are grateful to Prof. Ernst Hairer for helpful discussions related to
the example in Section \ref{sec:nlo} and for pointing out the closed formula
for $g(z)$ in \eqref{pofz}.

The first author was supported by the
King Abdullah University of Science and Technology (KAUST).
The second author was funded by the Deutsche Forschungsgemeinschaft (DFG, German Research Foundation)
under Germany's Excellence Strategy EXC 2044-390685587, Mathematics Münster:
Dynamics-Geometry-Structure.

\bibliographystyle{ACM-Reference-Format}
\bibliography{refs}

\end{document}